\documentclass[11pt]{amsart}
 \newtheorem{thm}{Theorem}[section]
 \newtheorem{cor}[thm]{Corollary}
 \newtheorem{lem}[thm]{Lemma}
 \newtheorem{prop}[thm]{Proposition}
 \theoremstyle{definition}
 \newtheorem{defn}[thm]{Definition}
 \newtheorem{rem}[thm]{Remark}
 \numberwithin{equation}{section}

\DeclareMathOperator{\sign}{sign}
\DeclareMathOperator{\Trace}{Trace}
\DeclareMathOperator{\Sign}{Sign}
\DeclareMathOperator{\Def}{def}
\begin{document}

\title[The Witten-Reshetikhin-Turaev Invariants of Lens Spaces]
{The Witten-Reshetikhin-Turaev Invariants of Lens Spaces}
\author{Khaled Qazaqzeh}

\address{Department of Mathematics, Yarmouk University, Irbid, 21163,
Jordan} \email{qazaqzeh@math.lsu.edu}
\urladdr{www.math.lsu.edu/\textasciitilde qazaqzeh/}
\date{08/29/2006}
\begin{abstract}
We derive an explicit formula for the Witten-Reshetikhin-Turaev $SO(3)$-invariants of lens spaces. We use the representation of the mapping class group of the torus corresponding to the Witten-Reshetikhin-Turaev $SO(3)$-TQFT to give such formula.
\end{abstract}
\maketitle
\section*{Introduction}
 We consider a variation of the $(2 + 1)$ cobordsim category that was studied in \cite{BHMV}. The variation consists of replacing the $p_{1}$ structure by integers that are called weights. This notion was first introduced by \cite{W,T}. This weighted category 
can be described roughly as follows.  The objects of this category $C$ are closed surfaces $\Sigma$ equipped with a Lagrangian subspace $\lambda  \subset H_{1}(\Sigma, \mathbb{R})$. We will denote objects by pairs $(\Sigma,  \lambda)$. A cobordism  from $(\Sigma , \lambda) $ to $(\Sigma', \lambda')$ is a 3-manifold 
with an orientation preserving homeomorphism (called its boundary identification) from its boundary to $-\Sigma\sqcup\Sigma'$. Here and elsewhere,   $-\Sigma$ denotes $\Sigma$ with the opposite orientation. Two cobordisms are equivalent if there is an orientation preserving homeomorphism between the underlying 3-manifolds that commutes with the boundary identifications.  A morphism  $M: (\Sigma , \lambda)  \rightarrow (\Sigma', \lambda')$ is  an equivalence class of cobordisms  from $(\Sigma , \lambda) $ to  $(\Sigma', \lambda')$ together with an integer weight.  We denote morphisms by $(M, w(M))$, where $w(M)$ denotes the weight of $M$. 
The gluing of 3-manifolds represents the composition and the weight of the composed morphism is given by ~\cite[Equation\,(1.6)]{GQ} which was derived from \cite[Thm.\,(4.1.1)]{T}.
 
The version of the WRT invariant that we study here is the invariant that is obtained from the $SO(3)$-TQFT-functor for $r\equiv 1 \pmod 4$ on $C$ over a commutative ring $K$ that was given in \cite{BHMV}. 
The ring is $K = \mathbb{Q}[\zeta_{2r},i]$ with $A=\zeta_{2r}$, and $\kappa=iA^{-1}$ where $\zeta_{2r} = e^{\frac{2\pi i}{2r}}$. 
This invariant is recovered from the representation of the gluing map between the two parts of the Heegaard splitting of the 3-manifold on the $SO(3)$-TQFT-vector space $V_{r}(\Sigma)$ of the boundary surface. 

I would like to acknowledge that the idea behind this work was inspired by Jeffery's work \cite{J}. The main difference between the two is the TQFT that was considered. In particular, her work was based on the TQFT associated to the group $SU(2)$ where our work is based on the TQFT associated to the group $SO(3)$. Also, I would like to thank my Ph.D adviser P. Gilmer for motivating and encouraging me to do this work.
\section{The $SO(3)$-TQFT-Representation of $SL(2,\mathbb{Z})$}

The $SO(3)$-TQFT associates to a surface $\Sigma$, a representation of the mapping class group on the vector space $V_{r}(\Sigma)$. We consider the case where the surface is the torus. The mapping class group of the torus is known to be $SL(2,\mathbb{Z})$. Hence, we obtain a representation of the group $SL(2,\mathbb{Z})$ on the vector space $V_{r}(S^{1}\times S^{1})$. The group $SL(2,\mathbb{Z})$ is specified in terms of two generators $\textbf{S} = \left(%
\begin{array}{cc}
  0 & -1 \\
  1& 0 \\
\end{array}%
\right)$, and 
$\textbf{T} = \left(%
\begin{array}{cc}
  1 & 1 \\
  0& 1 \\
\end{array}%
\right)$
with the relations $\textbf{S}^4 = (\textbf{S}\textbf{T})^6 = 1$. So any representation is specified in terms of these two generators satisfying the relations. 

We adopt the following notation of \cite{KM1} that has been used also in \cite{J}:
\[
e(\alpha)\stackrel{\Def}{=}\exp(2\pi i \alpha),
\]
\[
e_{n}(\alpha)\stackrel{\Def}{=}\exp(2\pi i \frac{\alpha}{n}).
\]
Also,
\[
\zeta\stackrel{\Def}{=}\exp\frac{i\pi}{4}.
\]

The representation is given explicitly by the assignments $\textbf{S}\mapsto \mathcal{D}^{-1} S$ and $\textbf{T}\mapsto T^{-1}$ as given in \cite[Page 98]{T} . This assignments define a projective representation because of the factor $\Delta \mathcal{D}^{-1}\in K$. This factor is expressed in terms of $\kappa$ in the following lemma.
\begin{lem}
Let $\kappa\in K$ be as above, then $\mathcal{D}\Delta^{-1} = \kappa^{3} = -iA^{-3}$.
\end{lem}
\begin{proof}
\begin{align*}
\kappa^{3}=\kappa^{3}\mathcal{D}\left\langle  S^{3},0\right\rangle &=\mathcal{D}\left(\kappa^{3}\left\langle S^{3},0\right\rangle\right)\\
&=\mathcal{D}\left\langle S^{3},1\right\rangle\\
&=\mathcal{D}\left(\mathcal{D}\Delta^{-1}\right)\left\langle S^{3},0\right\rangle\\
&=\mathcal{D}\left(\mathcal{D}\Delta^{-1}\right)\mathcal{D}^{-1}=\mathcal{D}\Delta^{-1}.
\end{align*}
\end{proof}
Hence the assignments $\textbf{S}\mapsto \mathcal{D}^{-1} S$ and $\textbf{T}\mapsto \kappa^{-1} T^{-1}$ give a linear representation as stated in the following proposition.
\begin{prop}
The representation of $SL(2,\mathbb{Z})$ on the vector space $V_{r}(S^{1}\times S^{1})$ obtained from the $SO(3)$-TQFT-representation is given by
\begin{align*}
S_{jl}&=\frac{1}{i\sqrt{r}}[e_{r}(jl)-e_{r}(-jl)]=\frac{2}{\sqrt{r}}sin(\frac{jl\pi}{r}),\\
T_{jl}&=\delta_{jl}T_{j},  \ \ \ T_{j}= ie_{2r}(j^{2}),
\end{align*}
for $1\leq j, l \leq \frac{r-1}{2}$. 
\end{prop}
\begin{rem}
In the above two matrices, we should have $j,l\in \{0,2,4,\ldots,{r-3}\}$ since the last set represents the set of simple objects in the $SO(3)-$modular category. However, the above two matrices and the actual matrices are similar by the same permutation matrix, i.e
\[
S = PS_{a}P^{-1}, \ \
T = PT_{a}P^{-1}
\] where $S_{a},T_{a}$ are the actual matrices.
\end{rem}
\begin{lem}\label{l:symmetry}
The coefficients $S_{jl},$ and $T_{j}$ satisfy the following symmetries:
\[
S_{jl}=S_{j(r+l)}=-S_{j(r-l)}, \ \ \ T_{j}=T_{(j+r)}=T_{(r-j)}.
\]
\end{lem}
\begin{proof}
It is clear for the coefficients in the $S$-matrix. For the coefficients in the $T$-matrix, we know that $e_{2r}(1)=-e_{r}(n)$ for some integer $n$. Therefore, $T_{j}=-ie_{r}(nj^2)=-ie_{r}(n(j+r)^2)=-ie_{r}(n(r-j)^2)$.
\end{proof}
Now, we wish to give an explicit formula for the representation of any element of $SL(2,\mathbb{Z})$ independent from the way we write that element in terms of the generators $\textbf{S}$, and $\textbf{T}$. But before we do that, we would like to quote the Gauss sum reciprocity formula in one dimension from \cite{J}.
\begin{prop}
If $\lambda, n, m \in \mathbb{Z}$ with $nm$ is even and $n\psi\in \mathbb{Z}$, then:
\begin{equation}\label{reciprocity}
\sum_{\lambda \pmod n} e_{2n}(m\lambda^{2})e(\psi \lambda) = \sqrt{\frac{in}{m}}\sum_{\lambda \pmod m} e_{2m}(-n(\lambda+\psi)^{2})
\end{equation}
\end{prop}

\subsection{The formula of the Representation}
Let $U = \left(%
\begin{array}{cc}
  a & b \\
  c& d \\
\end{array}%
\right) \in SL(2,\mathbb{Z}),$
we want to give a formula for the representation of $U$ in terms of its entries. 
\begin{defn}
A continued fraction expansion for $U$ is a tuple of integers $\mathcal{C}=(m_{1},\ldots,m_{t})$ such that
\[
U=\textbf{T}^{m_{t}}\textbf{S}\textbf{T}^{m_{t-1}}\ldots \textbf{T}^{m_{1}}\textbf{S}.
\]
\end{defn}
We would like to quote the following proposition from \cite{J}. 
\begin{prop}\label{jprop}\cite[Prop.\,(2.5)]{J}
Suppose $U$, and $\mathcal{C}$ as above. Then

$(i)\hspace{70pt} a/c = m_{t}-\cfrac{1}{m_{t-1} - \cfrac{1}{\ldots - \cfrac{1}{m_{1}}}}.$

$(ii) \hspace{60pt} b/a=-\left(\frac{1}{a_{1}} + \frac{1}{a_{2}a_{1}} + \ldots + \frac{1}{a_{t}a_{t-1}} \right).$

\noindent
Moreover, define $a_{i},b_{i},c_{i},d_{i}$ by the partial evaluation of this product:

\[\left(%
\begin{array}{cc}
  a_{i} & b_{i} \\
  c_{i}& d_{i} \\
\end{array}%
\right)=\textbf{T}^{m_{i}}\textbf{S}\textbf{T}^{m_{t-1}}\ldots \textbf{T}^{m_{1}}\textbf{S},\]

\noindent
with the convention that:
\[
a_{0}=d_{0}=1, \ \ b_{0}=c_{0}=0.
\]
Then these satisfy the recurrence relations (for $t\geq 2$)

$(iii) \hspace{40pt} a_{t} = m_{t}a_{t-1} - c_{t-1}, \hspace{30pt} c_{t} = a_{t-1};$

$(iv) \hspace{40pt} b_{t} = m_{t}a_{t-1} - d_{t-1}, \hspace{30pt} d_{t} = b_{t-1}.$

\end{prop}
\begin{lem}
Denote by $\Im_{t}$ the sum
\begin{equation}\label{e:sum}
\Im_{t}=\sum_{j_{1},\ldots,j_{t}=1}^{\frac{r-1}{2}} S_{j_{t+1}j_{t}}T_{j_{t}}^{m_{t}}S_{j_{t}j_{t-1}}T_{j_{t-1}}^{m_{t-1}}\ldots T_{j_{1}}^{m_{1}}S_{j_{1}j_{0}}.
\end{equation}
Then in terms of the previous proposition
\[
\Im_{t}= C_{t}\sum_{\stackrel{\gamma\pmod{2ra_{t}}}{\gamma \equiv j_{t+1}\pmod r}}\left\{e_{2ra_{t}}\left(-c_{t}(\gamma+\frac{j_{0}}{c_{t}})^{2}\right)-e_{2ra_{t}}\left(-c_{t}(\gamma-\frac{j_{0}}{c_{t}})^{2}\right)\right\}
,\]where \[
C_{t}=-i^{(t-1)}\zeta^{t-1}\zeta^{\sign(a_{t})}\frac{i^{(m_{t} + \dots + m_{1})}}{\sqrt{r\left|a_{t}\right|}}e_{2r}\left\{-\left(\frac{1}{a_{0}a_{1}}+\ldots+\frac{1}{a_{t-2}a_{t-1}}\right)j_{0}^2\right\}.
\]
\end{lem}
\begin{proof}
We observe that each of the indices $j_{1},\ldots,j_{t}$ appears in two of the $S$ coefficients in (\ref{e:sum}), so we may divide by 2 and replace the sum over $j_{1},\ldots,j_{t}$ from 1 to $r$ (using the symmetries in Lemma (\ref{l:symmetry})). To prove the result, we use induction on $t$.
For $t=1$:
\begin{align*}
\Im_{1}&=\frac{1}{2}\sum_{j_{1}=1}^{r}S_{j_{2}j_{1}}T_{j_{1}}^{m_{1}}S_{j_{1}j_{0}}\\
&=\frac{1}{2}\frac{-1}{r}i^{m_{1}}\sum_{j_{1}=1}^{r}e_{2r}(m_{1}j_{1}^2)\left\{e_{r}(j_{2}j_{1})-e_{r}(-j_{2}j_{1})\right\}\left\{e_{r}(j_{1}j_{0})-e_{r}(-j_{1}j_{0})\right\}\\
&=\frac{-1}{2r}2i^{m_{1}}\sum_{j\pmod r} e_{4r}(2m_{1}j^{2})\left\{e_{r}\left(\left(j_{2}+j_{0}\right)j\right)-e_{r}\left((j_{2}-j_{0})j\right)\right\}\\
&=\frac{-1}{r}i^{m_{1}}\sqrt{\frac{2ir}{2m_{1}}}\\
&\times \sum_{\beta\pmod {2m_{1}}}\left\{e_{4m_{1}}\left(-2r(\beta+\frac{j_{2}+j_{0}}{r})^{2}\right)-e_{4m_{1}}\left(-2r(\beta+\frac{j_{2}-j_{0}}{r})^{2}\right)\right\}\\
&=-i^{m_{1}}\zeta^{\sign(a_{1})}\frac{1}{\sqrt{r\left|a_{1}\right|}}\\
& \times \sum_{\beta\pmod {2a_{1}}}e_{2ra_{1}}\left(-\left(r\beta+j_{2}+j_{0}\right)^{2}\right)-e_{2ra_{1}}\left(-\left(r\beta+j_{2}-j_{0}\right)^{2}\right)
\end{align*}
Two terms in the third equality corresponding to the complex conjugates of the terms shown have been removed, and the overall expression was multiplied by 2: this results from the substitution $j\rightarrow -j$. The fourth equality was obtained by applying the reciprocity formula (\ref{reciprocity}). We obtain the required result  by substituting $\gamma=r\beta+j_{2}$. Hence, this confirms the first step of the induction.
Now, we assume that the result holds inductively for $t-1$. To prove the result for $\Im_{t}$, we use the symmetries in Lemma (\ref{l:symmetry}) to expand the sum over $j_{t}$ from 1 to $r$.
\begin{align*}
\Im_{t}&=\frac{1}{2}\sum_{j_{t}=1}^{r}S_{j_{t+1}j_{t}}T_{j_{t}}^{m_{t}}\Im_{t-1}\\
&=\frac{1}{2}\frac{i^{m_{t}}}{i\sqrt{r}}C_{t-1}\sum_{j_{t}\pmod {r}}\sum_{\stackrel{\gamma \pmod {2ra_{t-1}}}{\gamma \equiv j_{t} \pmod r}}e_{2ra_{t-1}}(m_{t}a_{t-1}j_{t}^{2}) e_{2ra_{t-1}}(-c_{t-1}\gamma^{2})\\& e_{2ra_{t-1}c_{t-1}}(-j_{0}^{2})\left\{e_{ra_{t-1}}(-\gamma j_{0})-e_{ra_{t-1}}(\gamma j_{0})\right\}\left\{e_{r}(j_{t+1}j_{t})-e_{r}(-j_{t+1}j_{t})\right\}\\
\end{align*}
Now, we replace $j_{t}$ by $\gamma$ and we combine the coefficients of the $\gamma^{2}$ factors using Proposition \ref{jprop}(iii), so we get:
\begin{align*}
\Im_{t}&=\frac{i^{m_{t}}}{2i\sqrt{r}}C_{t-1}e_{2ra_{t-1}c_{t-1}}(-j_{0}^{2})\sum_{\gamma\pmod{2rc_{t}}}e_{2ra_{t-1}}(a_{t}\gamma^{2}) \\
&\times \left\{e_{r}(j_{t+1}\gamma)-e_{r}(-j_{t+1}\gamma)\right\} \left\{e_{ra_{t-1}}(-\gamma j_{0})-e_{ra_{t-1}}(\gamma j_{0})\right\}.\\
\end{align*}
Now the substitution $\gamma \rightarrow -\gamma$ allows us to condense four terms to two and a factor of -2 in front of the sum, so this yields:
\begin{align*}
\Im_{t}&=\frac{-i^{m_{t}}}{i\sqrt{r}}C_{t-1}e_{2ra_{t-1}c_{t-1}}(-j_{0}^{2}) \sum_{\gamma\pmod{2rc_{t}}}e_{2rc_{t}}(a_{t}\gamma^{2})\\
& \times \left\{e_{rc_{t}}(c_{t}j_{t+1}+j_{0})\gamma)-e_{rc_{t}}(c_{t}j_{t+1}-j_{0})\gamma)\right\}\\
&=\frac{-i^{m_{t}}}{i\sqrt{r}}C_{t-1}e_{2ra_{t-1}c_{t-1}}(-j_{0}^{2})\sum_{\gamma\pmod{2rc_{t}}}e_{4rc_{t}}(2a_{t}\gamma^{2})\\
& \times \left\{e_{rc_{t}}(c_{t}j_{t+1}+j_{0})\gamma)-e_{rc_{t}}(c_{t}j_{t+1}-j_{0})\gamma)\right\}.\\
\end{align*}
Now, we use the reciprocity formula (\ref{reciprocity}) to obtain:
\begin{align*}
\Im_{t}&=\frac{-i^{m_{t}}}{i\sqrt{r}}C_{t-1}e_{2ra_{t-1}c_{t-1}}(-j_{0}^{2})\sqrt{\frac{2irc_{t}}{2a_{t}}}\sum_{\beta\pmod{2a_{t}}}\\
&\left\{e_{4a_{t}}\left(-2rc_{t}\left(\beta+\frac{c_{t}j_{t+1}+j_{0}}{rc_{t}}\right)^{2}\right)-e_{4a_{t}}\left(-2rc_{t}\left(\beta+\frac{c_{t}j_{t+1}-j_{0}}{rc_{t}}\right)^{2}\right)\right\}\\
&=\frac{-i^{m_{t}}}{i\sqrt{r}}C_{t-1}e_{2ra_{t-1}c_{t-1}}(-j_{0}^{2})\sqrt{\frac{irc_{t}}{a_{t}}}\sum_{\beta\pmod{2a_{t}}}\\
&\left\{e_{2ra_{t}}\left(-c_{t}\left(r\beta+j_{t+1}+\frac{j_{0}}{c_{t}}\right)^{2}\right)-e_{2ra_{t}}\left(-c_{t}\left(r\beta+j_{t+1}+\frac{-j_{0}}{c_{t}}\right)^{2}\right)\right\}.\\
\end{align*}
Now, the main formula holds using the substitution $\gamma=r\beta+j_{t+1}$. Finally, we use induction for $C_{t}$:
\[
C_{t}=\frac{-i^{m_{t}}}{i\sqrt{r}}C_{t-1}e_{2ra_{t-1}c_{t-1}}(-j_{0}^{2})\sqrt{\frac{irc_{t}}{a_{t}}}
=ii^{m_{t}}e_{2ra_{t-1}a_{t-2}}(-j_{0}^{2})\sqrt{\frac{ic_{t}}{a_{t}}}C_{t-1}.
\]
Hence, the formula for $C_{t}$ holds from the induction hypothesis for $C_{t-1}$.
\end{proof}

\begin{prop}
The representation of $SL(2,\mathbb{Z})$ on $V_{r}(S^{1}\times S^{1})$ is given by
\[
\mathcal{R}(U)_{jl}=(-iK_{t})\frac{1}{\sqrt{r\left|c\right|}}e_{2rc}(dl^{2})\sum_{\stackrel{\gamma \pmod {2rc}}{\gamma \equiv j_{t} \pmod r}}e_{2rc}(a\gamma^{2})\left\{e_{rc}(\gamma l)-e_{rc}(-\gamma l)\right\},
\]
where $K_{t} = i^{(t-1)}\zeta^{t-2}\zeta^{\sign(a_{t-1})}i^{(m_{t} + \dots + m_{1})}$ for $t\geq 2$ and
$K_{1}=i^{m_{1}}$.
\end{prop}
\begin{proof}
We prove the case for $t\geq 2$, we have
\begin{align*}
\mathcal{R}(U)_{j_{t}j_{0}}&=T_{j_{t}}^{m_{t}}\Im_{t-1}\\
&=i^{m_{t}}C_{t-1}\sum_{\stackrel{\gamma \pmod {2rc_{t}}}{\gamma \equiv j_{t} \pmod r}}e_{2r}(m_{t}\gamma^{2})\\
&\times \left\{e_{2ra_{t-1}}\left(-c_{t-1}(\gamma+\frac{j_{0}}{c_{t-1}})^{2}\right)-e_{2ra_{t-1}}\left(-c_{t-1}(\gamma-\frac{j_{0}}{c_{t-1}})^{2}\right)\right\}\\
&=i^{m_{t}}C_{t-1}e_{2ra_{t-1}c_{t-1}}(-j_{0}^{2})\sum_{\stackrel{\gamma \pmod {2rc_{t}}}{\gamma \equiv j_{t} \pmod r}}e_{2ra_{t-1}}(m_{t}a_{t-1}\gamma^{2})\\
&\times e_{2ra_{t-1}}(-c_{t-1}\gamma^{2})\left\{e_{ra_{t-1}}(-\gamma j_{0})-e_{ra_{t-1}}(\gamma j_{0})\right\}.\\
\end{align*}
\noindent
We combine the coefficients of the $\gamma^{2}$ factors using Proposition \ref{jprop}(iii), to obtain:
\begin{align*}
\mathcal{R}(U)_{j_{t}j_{0}}&=-i^{m_{t}}C_{t-1}e_{2ra_{t-1}a_{t-2}}(-j_{0}^{2})\sum_{\stackrel{\gamma \pmod {2rc_{t}}}{\gamma \equiv j_{t} \pmod r}}\\
&e_{2ra_{t-1}}(a_{t}\gamma^{2}) \left\{e_{ra_{t-1}}(\gamma j_{0})-e_{ra_{t-1}}(-\gamma j_{0})\right\}\\
&=-i^{m_{t}}C_{t-1}e_{2ra_{t-1}a_{t-2}}(-j_{0}^{2})e_{2ra_{t}a_{t-1}}(-j_{0}^{2})\sum_{\stackrel{\gamma \pmod {2rc_{t}}}{\gamma \equiv j_{t} \pmod r}}\\
&\left\{e_{2ra_{t-1}}\left(a_{t}\left(\gamma+\frac{j_{0}}{a_{t}}\right)^{2}\right)-e_{2ra_{t-1}}\left(a_{t}\left(\gamma-\frac{j_{0}}{a_{t}}\right)^{2}\right)\right\}.
\end{align*}
\noindent
We substitute the value of $C_{t-1}$ from the previous lemma, to get:
\begin{align*}
\mathcal{R}(U)_{j_{t}j_{0}}&=i^{m_{t}}i^{(t-2)}\zeta^{t-2}\zeta^{\sign(a_{t-1})}\frac{i^{(m_{t-1} + \dots + m_{1})}}{\sqrt{r\left|a_{t-1}\right|}}\\
&e_{2r}\left\{-\left(\frac{1}{a_{0}a_{1}}+\ldots+\frac{1}{a_{t-3}a_{t-2}}\right)j_{0}^2\right\}
 e_{2ra_{t-1}a_{t-2}}(-j_{0})e_{2ra_{t}a_{t-1}}(-j_{0}^{2})\\
&\times \sum_{\stackrel{\gamma \pmod {2rc_{t}}}{\gamma \equiv j_{t} \pmod r}}\left\{e_{2ra_{t-1}}\left(a_{t}\left(\gamma+\frac{j_{0}}{a_{t}}\right)^{2}\right)-e_{2ra_{t-1}}\left(a_{t}\left(\gamma-\frac{j_{0}}{a_{t}}\right)^{2}\right)\right\}.\\
\end{align*}
Now, we use Proposition \ref{jprop}(ii)for the prefactor involving $j_{0}^{2}$ to obtain:
\begin{align}\label{rep}
\mathcal{R}(U)_{j_{t}j_{0}}&=-iK_{t}\frac{1}{\sqrt{r\left|c\right|}}e_{2ra}(bj_{0}^{2})\\
&\sum_{\stackrel{\gamma \pmod {2rc_{t}}}{\gamma \equiv j_{t} \pmod r}}\left\{e_{2ra_{t-1}}\left(a_{t}\left(\gamma+\frac{j_{0}}{a_{t}}\right)^{2}\right)-e_{2ra_{t-1}}\left(a_{t}\left(\gamma-\frac{j_{0}}{a_{t}}\right)^{2}\right)\right\}\notag
\end{align}
\begin{align*}
&=-iK_{t}\frac{1}{\sqrt{r\left|c\right|}}e_{2rac}(bcj_{0}^{2})e_{2rac}(j_{0}^{2}) \sum_{\stackrel{\gamma \pmod {2rc_{t}}}{\gamma \equiv j_{t} \pmod r}}e_{2rc}(a\gamma^{2})\left\{e_{rc}(\gamma j_{0})-e_{rc}(-\gamma j_{0})\right\}\notag\\
&=-iK_{t}\frac{1}{\sqrt{r\left|c\right|}}e_{2rc}(dj_{0}^{2})\sum_{\stackrel{\gamma \pmod {2rc}}{\gamma \equiv j_{t} \pmod r}}e_{2rc}(a\gamma^{2})\left\{e_{rc}(\gamma j_{0})-e_{rc}(-\gamma j_{0})\right\}.
\end{align*}
Here $K_{t}=i^{(t-1)}\zeta^{t-2}\zeta^{\sign(a_{t-1})}i^{(m_{t} + \dots + m_{1})}$. The last equality follows from the equation $bc+1=ad$.
\end{proof}

\section{The formula of the invariant for Lens Spaces}
The lens space $L(p,q)$ is specified by a pair of coprime integers $p,q$. We can assume that $0<-q<p$ as it was shown in \cite{R} that $L(p,q)$ is diffeomorphic to $L(p,q+np)$ for any integer $n$. The above lens space is obtained by doing a rational surgery on $S^{3}$ along the unknot with coefficient $-p/q$. Equivalently, it is obtained by an integer surgery on $S^{3}$ along a chain link $L$ with successive framings by integers $m_{1}, m_{2}, \ldots, m_{t-1}$ such that $\mathcal{C}=( m_{1}, m_{2}, \ldots, m_{t-1})$ is a continued fraction of $-p/q$ as in \cite{PS}. The rational surgery on $S^{3}$ means removing a solid torus around the unknot and gluing it back using the matrix:
\[
A = \left(%
\begin{array}{cc}
  p & d \\
  -q& -b \\
\end{array}%
\right) \in SL(2,\mathbb{Z}).\]

\noindent
Hence, the lens space $L(p,q)$ is obtained by gluing two solid tori by

\[ U= \textbf{S}A = \left(%
\begin{array}{cc}
  q & b \\
  p& d \\
\end{array}%
\right).\]

Therefore, $U$ has a continued fraction expansion $\mathcal{C} = (m_{1}, \ldots, m_{t-1}, m_{t} =0 )$. Thus, according to the TQFT axioms, the WRT invariant of $L(p,q)$ with the weight  obtained from the composition  is given by 
\[\left\langle L(p,q),w\right\rangle_{r} = \mathcal{R}(U)_{11},\]
where $w = \sum_{i=2}^{t} \sign{(c_{i-1}c_{i})}$ as computed in \cite[Page.\,415]{G}.

As $0<-q<p$, we have $-p/q >1$. We need to use a generalized version of \cite[Lem.\,(3.1)]{J}.
\begin{lem}
$-p/q$ has a unique continued fraction expansion with all $m_{i}\geq 2$.
\end{lem} 
\begin{proof}
We know that any rational number has a continued fraction expansion. The construction of the required fraction expansion is given as follows: set $m_{t-1}=\left\lceil -p/q \right\rceil\geq 2$, so $-p/q = m_{t-1} - (p^{'}/q^{'})^{-1}$, where $1 \leq q^{'} < p^{'}= -q < p$, and $p^{'}/q^{'} > 1$. We repeat this process for $p^{'}/q^{'}$ and since the denominators continue to decrease, this process will terminate. This gives a continued fraction expansion with $m_{i}\geq 2$. The idea of this construction is due to Jeffery in the her proof of \cite[Lem.\,(3.1)]{J}. To prove the uniqueness, we assume that $\mathcal{C}^{'} = (m_{1}, m_{2}, \ldots, m_{r-1})$ is another continued fraction expansion of $-p/q$ with all $m_{i}\geq 2$. Let $k$ be the first integer such that $m_{t-1-k} \neq m_{r-1-k}$. Hence, we would have $s/u = m_{r-1-k} - (s^{'}/u^{'})^{-1}$ with $s^{'}/u^{'} < 1$. Therefore, $(m_{1}, m_{2}, \ldots, m_{r-2-k})$ is a continued fraction expansion for $ (s^{'}/u^{'}) < 1$ with $m_{i} \geq 2$ which contradicts the result of the next lemma.
\end{proof}
\begin{lem}
If $s/u < 1$, then there is no continued fraction expansion for $s/u$ with all $m_{i} \geq 2$.
\end{lem}
\begin{proof}
Assume such a continued fraction expansion exists for $s/u$. We consider first the case $0 < s/u < 1$. If $s/u = m_{v-1} - (s^{'}/u^{'})^{-1}$ with $ 1 \leq s < u = s^{'} < u^{'}$ then it implies $0 < s^{'}/u^{'} < 1$. Now, if we repeat this process for $s^{'}/u^{'}$
then it will never terminate as long as $m_{i} \geq 2$. For the second case, if $s/u < 0$ then it is enough to notice that $ s/u = m_{v-1} - (s^{'}/u^{'})^{-1}$ where $ 0 < s^{'}/u^{'} < 1$. Now we can apply the first case to $s^{'}/u^{'}$ to conclude that there is no continued fraction expansion with all $m_{i} \geq 2$.
\end{proof}
\begin{cor}
There is a unique continued fraction expansion $\mathcal{C}$ for any rational number $\left| s/u \right| \neq 1$ with 
\begin{enumerate}
	\item $m_{i} \geq 2, \   1\leq i \leq t$ if $s/u > 1$.
	\item $m_{i} \leq -2, \   1\leq i \leq t$ if $s/u < -1$.
	\item $m_{i} \leq -2, \   1\leq i \leq t-1$ and $m_{t} = 0$ if  \ $0 \leq s/u < 1$.
	\item $m_{i} \geq 2, \   1\leq i \leq t-1$ and $m_{t} = 0$ if  \ $-1 < s/u \leq 0$.
\end{enumerate}
\end{cor}
\noindent
We need to use the next lemma whose proof can be found in \cite{J}.
\begin{lem}\cite[Lem.\,(3.2)]{J}
If $A$ is given by the continued fraction with all $m_{i} \geq 2$ as above, then
\[
\Phi(U) = 
-3(t-1) + \sum_{i = 0}^{t-1} m_{i}.
\]

\end{lem}
Here and elsewhere, $\Phi(U)$ is the Rademacher phi function of $U$ (see \cite{J,KM2} for more details about this function).
Now, we have
\[
\kappa^{3\Sign(W_{L})}\left\langle L(p,q),0\right\rangle_{r}= \left\langle L(p,q),w\right\rangle_{r} =\kappa^{\Trace(W_{L})}\mathcal{R}(U)_{11}.\\
\]
where $W_{L}$ is the linking matrix of the link $L$.
Hence, we have
\begin{align}\label{formula}
\left\langle L(p,q),0\right\rangle_{r}=\kappa^{\Phi(U)}\mathcal{R}(U)_{11}.
\end{align}
Therefore, we can conclude the following theorem.
\begin{thm}
The Witten-Reshetikhin-Turaev invariant of the lens space $L(p,q)$ weighted zero is given by
\[
\left\langle L(p,q),0\right\rangle_{r} = -\frac{i\zeta^{t-1}}{\sqrt{rp}}e_{2r}(12s(q,p))\sum_{\pm}\sum_{n=1}^{p}e_{rp}(\pm 1)e_{p}(2qrn^{2})e_{p}(2n(q\pm1)),
\]
where the Dedekind sum $s(q,p)$ is defined in \cite[Equation.\,(2.16)]{J}.
\end{thm}
\begin{proof}
Using equation (\ref{rep}), we have
\begin{align*}
\mathcal{R}(U)_{11}&=-iK_{t}\frac{1}{\sqrt{rp}}e_{2rq}(b)\\
&\sum_{\stackrel{\gamma \pmod {2rp}}{\gamma \equiv 1 \pmod r}}\left\{e_{2rp}\left(q\left(\gamma+\frac{1}{q}\right)^{2}\right)-e_{2rap}\left(q\left(\gamma-\frac{1}{q}\right)^{2}\right)\right\}\\
&=-iK_{t}\frac{1}{\sqrt{rp}}e_{2rq}(b)\sum_{n=1}^{p}\\ & \left\{e_{2rpq}(2rnq+q+1)^{2}-e_{2rpq}(2rnq+q+1)^{2}\right\}\\
&=-iK_{t}\frac{1}{\sqrt{rp}}e_{2rq}(b)\sum_{n=1}^{p}e_{p}(2qrn^{2})\\
& \left\{e_{p}(2n(q+1))e_{2rpq}(q+1)^{2}-e_{p}(2n(q-1))e_{2rpq}(q-1)^{2}\right\}.
\end{align*}
Therefore by (\ref{formula}), we get:
\begin{align*}
\left\langle L(p,q), 0\right\rangle_{r} &=-ii^{\Phi(U)}\zeta^{t-2}\zeta^{\sign(p)}\left(i^{-1}e_{2r}(-1)\right)^{\Phi(U)}\frac{1}{\sqrt{rp}}e_{2rq}(b)\\
& \times \sum_{n=1}^{p}e_{p}(2qrn^{2})\\&\times \left\{e_{p}(2n(q+1))e_{2rpq}(q+1)^{2}-e_{p}(2n(q-1))e_{2rpq}(q-1)^{2}\right\}.
\end{align*}
Now,
\[
e_{2rq}(b)e_{2rpq}(q\pm 1)^{2}=e_{2rpq}(bp + q^2 \pm 2q + 1)=e_{2rp}(d + q \pm 2),
\]as  $bp + 1 = dq$.
Also, if we introduce the integer $q^{*}$ solving $qq^{*} \equiv 1 \pmod{p}$ then we have
\[
e_{2rp}(d + q \pm 2)\left(e_{2r}(-1)\right)^{\Phi(U)}=e_{rp}(\pm 1)e_{2r}(12s(d,p))=e_{rp}(\pm 1)e_{2r}(12s(q^{*},p)),
\] as $\Phi(U) = \frac{d + q}{p} - 12 s(d,p)$ and $q^{*} \equiv d \pmod{p}$. Finally, we obtain the result as $s(q^{*},p) = s(q,p)$.
\end{proof}
\bibliographystyle{amsalpha}

\end{document}